\documentclass[letterpaper, 10 pt, conference]{ieeeconf}  % Comment this line out if you need a4paper

\IEEEoverridecommandlockouts                              % This command is only needed if 
\usepackage{graphicx}
\usepackage{epstopdf}
\usepackage{epsfig} % for postscript graphics files
\usepackage{mathptmx} % assumes new font selection scheme installed
\usepackage{amsmath} % assumes amsmath package installed
\usepackage{amssymb}  % assumes amsmath package installed
\usepackage{hyperref}
\usepackage{siunitx}

\title{\LARGE \bf
Optimal-power Configurations for Hover Solutions in Mono-spinners
}

\author{Mojtaba Hedayatpour$^{1}$, Mehran Mehrandezh$^{1}$ and Farrokh Janabi-Sharifi$^{2}$% <-this % stops a space
\thanks{$^{1}$Mojtaba Hedayatpour and Mehran Mehrandezh are with the Faculty of Engineering \& Applied Science, University of Regina, 3737 Wascana Parkway, Regina, Saskatchewan, Canada, S4S 0A2.
        {\tt\small \{hedayatm,mehran.mehrandezh\}@uregina.ca}}%
\thanks{$^{2}$Farrokh Janabi-Sharifi is with the Department of Mechanical \& Industrial Engineering, Ryerson University,
        350 Victoria Street, Toronto, Ontario, Canada, M5B 2K3.
        {\tt\small fsharifi@ryerson.ca}}%
}

\begin{document}

\maketitle
\thispagestyle{empty}
\pagestyle{empty}

%%%%%%%%%%%%%%%%%%%%%%%%%%%%%%%%%%%%%%%%%%%%%%%%%%%%%%%%%%%%%%%%%%%%%%%%%%%%%%%%
\begin{abstract}

Rotary-wing flying machines draw attention within the UAV community for their in-place hovering capability,
and recently, holonomic motion over fixed-wings. In this paper, we investigate about the power-optimality in a mono-spinner, i.e., a class of rotary-wing UAVs with one rotor only, whose main body has a streamlined shape for producing additional lift when counter-spinning the rotor. We provide a detailed dynamic model of our mono-spinner. Two configurations are studied: (1) a symmetric configuration, in which the rotor is aligned with the fuselage's COM, and (2) an asymmetric configuration, in which the rotor is located with an offset from the fuselage's COM. While the
former can generate an in-place hovering flight condition, the latter can achieve trajectory tracking in 3D space by resolving the yaw and precession rates. Furthermore, it is shown that by introducing a tilting angle between the rotor and the fuselage, within the asymmetric design, one can further minimize the power consumption without compromising the overall stability. It is shown that an energy optimal solution can be achieved through the proper aerodynamic design of the mono-spinner for the first time.

\end{abstract}

%%%%%%%%%%%%%%%%%%%%%%%%%%%%%%%%%%%%%%%%%%%%%%%%%%%%%%%%%%%%%%%%%%%%%%%%%%%%%%%%
\section{INTRODUCTION}

Unmanned Aerial Vehicles (UAVs) have gained a significant amount of attention within the industry and academia in recent years. Rotary-wing UAVs, in particular, are gaining attention within the research community for their leading advantage over fixed-wings for their in-place hovering capability, holonomy in motion, and safe operation due to their redundancy. However, they still fall behind in terms of the flight time and power consumption. They are used in a variety of applications such as: inspection of infrastructures, object delivery, sports broadcasting and even for calibrating the antenna of a radio telescope \cite{delivery}-\cite{telescope}. 

Under- and also fully-actuated rotary-wing UAVs have been cited in the literature, \cite{omnicopter}-\cite{underac-prop}. Some examples of the under-actuated hover-capable  rotary-wing UAVs are: quad/hexa-copters, co-axial helicopters, mono-spinners and ornithopters, \cite{underac-prop}-\cite{monospinner}. Under-actuated UAVs have become popular due to their minimalistic and simple design. They can be classified under three different categories: (1) Samara-type, (2) flapping wing, and (3) spinners.
The vehicles in the first category are inspired from nature through unpowered flight in maple/pine seeds (or Samaras), \cite{samara}, \cite{samaraCityU}. They offer a passively-stable flight with slow-rate descending altitude, therefore they would require no active control. A small propulsion system, however, has been added to the body, in lab-scale prototypes, to control the rotational speed and correspondingly the descending rate \cite{samara}.

The vehicles in the second category are inspired from birds species. They normally have one or two flapping wings. Single actuator (wing) type vehicles of this category are only capable of altitude control, while those with two actuators (i.e., two wings) are capable of controlling all three translational degrees of motion, \cite{flapwing}. 

Flying machines falling into the third category, namely spinners, generate a constant rotation about a fixed axis in space, i.e., precession axis, due to the presence of unbalanced moments in the system (aka, boomerang-type spinners). The minimum number of actuators required to achieve position control in these machines is one, which renders itself as the simplest structure of rotary-wing UAV that is controllable in all translational degrees of motion, \cite{spinners}. 

Technical specifications of a spinner-type machine with only one actuator, (aka, a mono-spinner) can be found in \cite{monospinner}, and also specifications for another spinner with two actuators can be found in \cite{spinners}. A novel design of a small spinning vehicle that consists of a single propeller and an aerodynamically-designed streamline-shape fuselage is presented in \cite{picco}. A single blade spinning rotor-craft with two tilted rotors is presented in \cite{phantom}. A transformable vehicle, also known as THOR, which transforms from a fixed-wing aircraft to a spinner can be found in \cite{thor}.

\begin{figure}[b]
\centering
\includegraphics[scale=0.25]{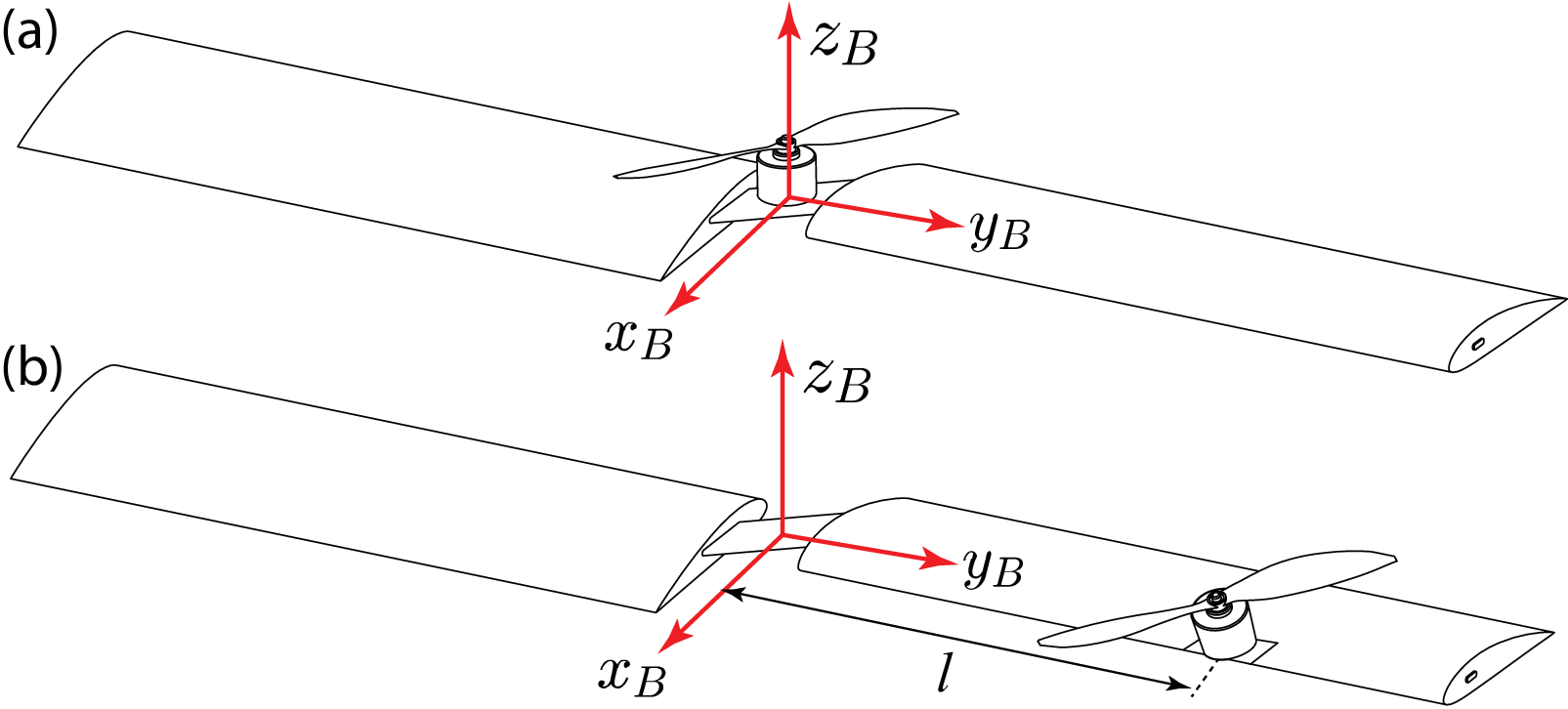}
\caption{ A spinning UAV with a single actuator and streamline-shape aerodynamic fuselage. a) The actuator is located at the center of mass (COM) of the vehicle. b) the actuator is located at distance \emph{l} from COM of the vehicle and is titled about the y-axis of the body frame. The red frame denotes the body frame and is attached to the COM of the vehicle. }
\label{myUAV}
\end{figure}

This paper investigates about the optimal design of mono-spinners with a streamline-shape  fuselage that can provide the most energy-efficient solution for hovering and also position control for the first time. For simplicity, in this paper, we only consider blade like geometry for the fuselage. However, any streamline-shape fuselage could be considered as long as the aerodynamic model for forces and moments is available and could be incorporated in the mathematical modeling of the UAV. 
 
A comprehensive dynamic model of a mono-spinner whose rotor is positioned at an offset from the fuselage's COM is developed considering the ``blade element theory", \cite{anderson}. This dynamic model is then used to find the optimal configuration, i.e., optimal aerodynamic characteristics of the rotor and the streamline-shape fuselage, for inducing the largest lift and smallest drag simultaneously. More specifically, two possible configurations are explored: (1) a symmetric configuration, in which the rotor is located at the fuselage's COM, and (2) an asymmetric configuration, in which the rotor is located with an offset w.r.t the fuselage's COM. Only, in-place hovering can be achieved under (1), however, position control is possible under (2). 

Optimal configurations that would yield the most efficient in-position hovering for both configurations are formulated and design guidelines are provided. The optimal solutions, however, were obtained for a steady-state flight (no position control). It is noteworthy that in the second configuration, a pseudo in-place hovering is achieved in which case yawing and precession speeds of the spinner will not be zero, but bounded. Furthermore, it is shown that power  consumption for the flight can be further reduced by introducing a tilting angle between the rotor and the fuselage's principle axis.

This paper is organized as follows: Section II presents mathematical modeling. Section III presents optimal hover solutions for two aforementioned configurations for mono-spinners. Section IV presents a discussion on results obtained via two possible configurations, while pinpointing  key characteristics of the two possible designs. Finally, conclusions and future works are presented in section V.

\begin{figure}[t]
\centering
\includegraphics[scale=0.5]{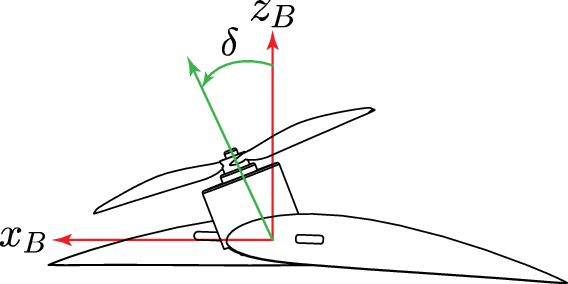}
\caption{ The rotor can be tilted about y-axis of the body frame by angle $\delta$. The positive direction of $\delta$ is shown in green. }
\label{tilt_angle}
\end{figure}

\section{MATHEMATICAL MODELING}

This section presents mathematical modeling of the proposed vehicle considering the effects of tilting in the rotor, having nonzero freestream velocity and a streamlined fuselage. Before getting into equations, some explanations and clarifications used via modeling are given below.  

\begin{figure}[b]
\centering
\includegraphics[scale=0.28]{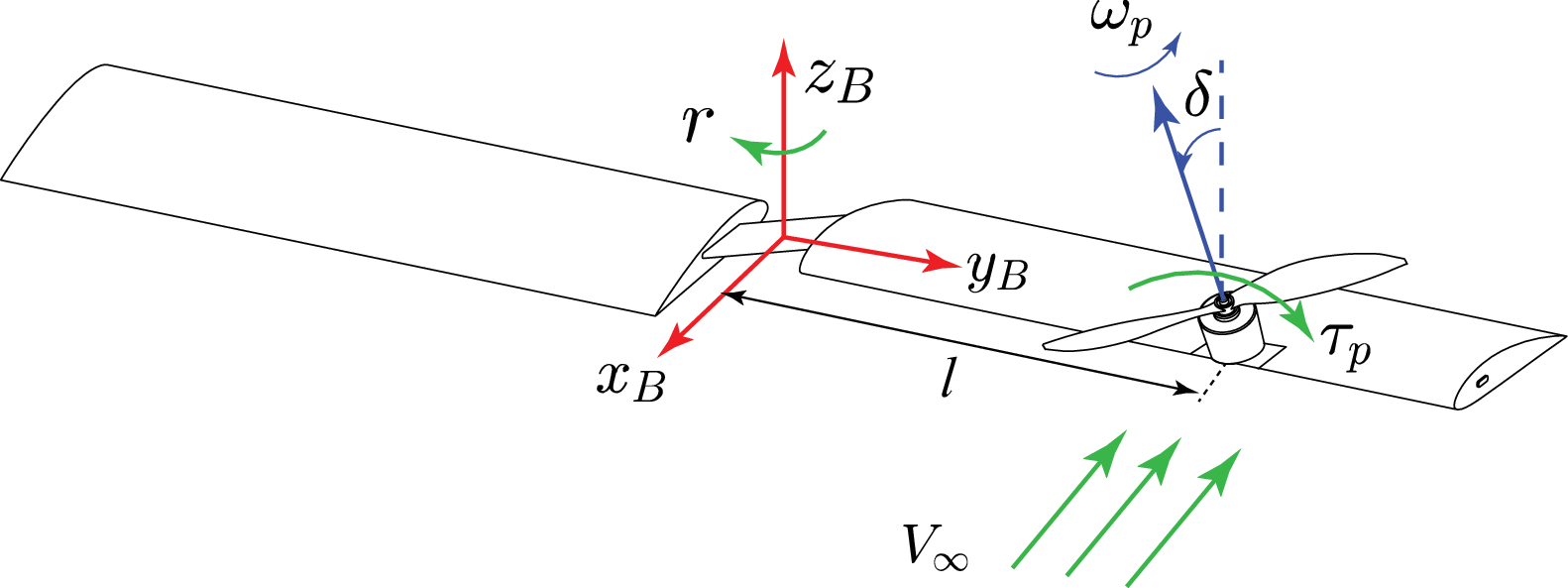} 
\caption{ The vehicle has a yaw rate $r$ and the propeller is located at distance $l$ from the COM of the vehicle and is turning with angular velocity $\omega_p$ as shown in blue. The propeller experiences an almost uniform freestream velocity $V_{infty}=rl$ (shown in green) which generates the moment $\tau_p$ due to asymmetrical lift distribution over the propeller. The propeller is also tilted by angle $\delta$. }
\label{Vinf_effect}
\end{figure}

Firstly, freestream velocity can have significant effects on propeller's performance. In spinners, when the propeller is located at a distance from COM of the vehicle and as the body is spinning, the COM of the propeller goes through a rotation about the COM of the vehicle and it experiences an almost uniform freestream velocity of magnitude $V_{\infty}=rl$ as shown in Fig. \ref{Vinf_effect}. Depending on orientation of the propeller, the freestream velocity may have two components: (1) one parallel to the axis of rotation of the propeller, and (2) one perpendicular to the axis of rotation of the propeller. According to blade element theory (BET) \cite{anderson}, the former component changes the angle of attack of the blade element with respect to the local airflow. If freestream velocity is constant, this change in angle of attack can be compensated by adding it to the original angle of attack of the blade element (when freestream velocity is zero). The latter component changes the magnitude of the local airflow velocity over each blade element, which could be either positive (for the advancing blade) or negative (for the retreating blade) and thus it generates asymmetrical lift distribution over the propeller and a rolling moment $\tau_p$ as shown in Fig. \ref{Vinf_effect}. Detailed analysis of these effects can be found in \cite{iros}-\cite{cdsr}. 

Secondly, designing an aerodynamic fuselage for the vehicle to maximize lift versus drag can help minimize the power consumption. In this paper, we assume that the fuselage is streamline-shape as shown in Fig. 1. The reaction moment of the propeller and cross-coupling effects of angular momentum in the system can spin the streamline-shape fuselage and generate extra lift which may be used to balance weight partially and thus saving some energy. 

Lastly, the effects of location and orientation of the propeller are incorporated into equations. The mathematical model is derived considering all of the above-mentioned effects. It is shown that by selecting appropriate values for the parameters of the system (e.g., tilting angle of the rotor, distance of the propeller from COM of the vehicle and etc.), mathematical model for some specific configurations (e.g., co-axial helicopter and monospinner) can be obtained. 

\subsection{Notation \& Parameters}

Throughout this paper, boldface letter \textbf{R} is used to represent rotation matrices. Two reference frames are introduced: (1) an inertial frame that is attached to the earth and is represented by $I$, and (2) a body frame that is attached to the center of mass of the vehicle and rotates with it, represented by $B$. The transformation matrix from body frame to the inertial frame is shown by $^I$\textbf{R}$_B$. The angular velocity of the vehicle, as expressed in body frame, is represented by $\pmb{\omega}_B=(p,q,r)$. The tilting angle of the rotor as shown in Fig. \ref{tilt_angle}, is represented by $\delta$ in radians. Greek letter $\alpha$ is reserved to represent angle of attack of an airfoil. 

\subsection{Equations of Motion}

The vehicle consists of several rigid bodies including the fuselage, and the propeller. The propeller is located at distance $l$ from the COM of the vehicle along the y-axis and is tilted by angle $\delta$. The propeller has two blades with chord $c_p$ and has a radius $R_p$. The fuselage is assumed to be streamline-shape, aerodynamic and similar to the propeller with two blades with chord $c_B$ and radius $R_B$. For simplicity and to avoid adding more complexity, we keep the airfoil parameters fixed (NACA 4415 is selected as the airfoil in both propeller and the fuselage~\cite{uiuc}). 

It is assumed that the vehicle is symmetric about its axes of rotation, therefore, its moment of inertia tensor will be diagonal. Moment of inertia tensors for the spinning fuselage and rotor are modeld as that in a rotating disk, and are represented by $I^p$ \& $I^B$ respectively. Also, the additional moment of inertia due to the battery and electronics, attached to the fuselage, is represented by $I^e$, and is modeled as that in a cylinder centered at the COM of the vehicle. Therefore, total moment of inertia of the vehicle can be calculated as follows: 
  \begin{equation}
      I^T = I^e+I^p+I^B
  \end{equation}

A DC motor with a torque constant $K_{\tau_m}$, electric resistance $R_m$, electric inductance $L_m$, electric current $i$, voltage $V_m$  and electromotive force constant $K_v$ turns the propeller. Dynamic equations of the motor can be written as follows: 
  \begin{equation}
  \label{motor_eq}
  \begin{matrix}
      I^p\dot{\omega_p}+\tau_{d_p}=K_{\tau_m}i \\
      L_m\frac{di}{dt}+R_mi=V_m-K_v\omega_p
  \end{matrix}
  \end{equation}

The power consumption of the motor can be calculated as follows: 
\begin{equation}
	P_m = V_mi
\end{equation}

The equations governing the rotational motion of the vehicle, as expressed in the body frame, can be written as follows: 
  \begin{equation}
  \label{rot_eq}
      I^T\pmb{\dot{\omega}}_B + I^p\pmb{\dot{\omega}}_p + sk(\pmb{\omega}_B)(I^T\pmb{\omega}_B + I^p(\pmb{\omega}_p+\pmb{\omega}_B)) = \pmb{\tau}_{ext}
  \end{equation}
where $sk(\pmb{\omega}_B)$ is the skew-symmetric matrix of the angular velocity vector of the vehicle. $\pmb{\tau}_{ext}$ represents total external moments applied to the vehicle as follows: 
  \begin{equation}
      \pmb{\tau}_{ext} = \pmb{\tau}_{f} + \pmb{\tau}_{p} + \pmb{\tau}_{d_p} + \pmb{\tau}_{d_B}
  \end{equation}
where $\pmb{\tau}_{f}, \pmb{\tau}_{p}, \pmb{\tau}_{d_p}$ and $\pmb{\tau}_{d_B}$ represent moments due to thrust force of the propeller, asymmetrical lift distribution over the propeller, drag force of the propeller and drag force of the fuselage respectively. According to \cite{cdsr} (details of the derived equations can be found in~\cite{cdsr}), using blade element theory, these moments can be calculated as follows: 
\begin{equation}
\begin{matrix}
	\pmb{\tau}_{f} = \pmb{f}_pl
\end{matrix}
\end{equation}
\begin{equation}
\begin{matrix}
    \pmb{f}_p = \rho c_p C_{L_p} \Big( \frac{R_p^3\omega_p^3}{3} + \frac{R_p^3r^2\cos^2\delta}{3} + \frac{R_pr^2l^2\cos^2\delta}{2} + \\
    \frac{2R_p^3r\omega_p\cos\delta}{3} \Big) \pmb{e}_{f_p}
\end{matrix}
\end{equation}
\begin{equation}
	\pmb{\tau}_{p} = \rho c_p C_{L_p} \Big( \frac{R_p^3\omega_prl\cos\delta + R_p^3r^2l\cos\delta}{3} \Big) \pmb{e}_{\tau_p}
\end{equation}
\begin{equation}
	\begin{matrix}
	\pmb{\tau}_{d_p} = \rho c_p C_{D_p} \Big( \frac{R_p^4\omega_p^3}{4} + \frac{R_p^4r^2\cos^2\delta}{4} + \frac{R_p^2r^2l^2\cos^2\delta}{2} + \\
    \frac{R_p^4r\omega_p\cos\delta}{2} \Big) \pmb{e}_{\tau_{dp}}
	\end{matrix}
\end{equation}
\begin{equation}
	\pmb{\tau}_{d_B} = \Big(\frac{1}{4}\rho c_B C_{D_B} R_B^4r^2-\gamma{r}\Big) \pmb{e}_{\tau_{dB}}
\end{equation}
where $l$ represents the distance of the COM of the propeller from the COM of the vehicle and $\gamma$ represents drag coefficient of the central hub of the vehicle which includes battery and electronics. $\pmb{e}_{f_p}, \pmb{e}_{\tau_p}, \pmb{e}_{\tau_{dp}}$ and $\pmb{e}_{\tau_{dB}}$ are appropriate unit vectors to determine direction of the forces and moments expressed in the body frame. $\rho$ is the air density, $c_p, c_B, R_p$ and $R_B$ represent chord and radius of the propeller and fuselage respectively. $C_{L_p}, C_{L_B}, C_{D_p}$ and $C_{D_B}$ represent aerodynamic coefficients of the airfoil used in propeller and the streamline-shape fuselage. NACA 4415 airfoil is used for both propeller and fuselage. These aerodynamic coefficients as a function of angle of attack can be obtained from experimental results from \cite{uiuc}. Assuming the angle of attack varies between zero and 10 degrees, the following linear functions can be written: 
\begin{equation}
\begin{matrix}
	C_{L_p} = 0.1\alpha_p + 0.5 \:\:,\:\: C_{L_B} = 0.1\alpha_B + 0.5 \\
    C_{D_p} = 0.006\alpha_p + 0.04 \:\:,\:\: C_{D_B} = 0.006\alpha_B + 0.04
\end{matrix}
\end{equation}
where $\alpha_p$ and $\alpha_B$ are the angles of attack for the propeller and the streamline-shape fuselage respectively. 

The position of the vehicle expressed in the inertial frame $I$ is denoted by $\pmb{d}$. Using Newton's second law, the equation governing translational motion of the vehicle can be written as follows: 
\begin{equation}
	m\ddot{\pmb{d}} = m\pmb{g} + |\pmb{f}_p|\pmb{e}_{f_p}^I + |\pmb{f}_B|\pmb{e}_{f_B}^I
\end{equation}
where $m$ is the total mass of the vehicle and $\pmb{g}$ is the gravitational acceleration as expressed in the inertial frame. $\pmb{e}_{f_p}^I$ and $\pmb{e}_{f_B}^I$ are the unit vectors expressed in the inertial frame determining the direction of thrust force of the propeller and the fuselage. Thrust force of the streamline-shape fuselage, $\pmb{f}_B$, can be expressed in the body frame as follows: 
\begin{equation}
	\pmb{f}_B = \frac{1}{3}\rho c_B C_{L_B} R_B^3r^2 \pmb{e}_{f_{B}}
\end{equation}
where $\pmb{e}_{f_{B}}$ is the unit vector determining the direction of $\pmb{f}_B$ in the body frame. Note that it is assumed the translational velocities are slow enough to neglect drag forces in translational motion. 

\subsection{Hover Solution}

In our spinner, hovering is defined as the state at which the vehicle is turning around a fixed axis in space with constant angular velocity and bounded linear velocities \cite{spinners}. 

Assume that the vehicle is hovering and rotating about a fixed axis $\pmb{n}$. The evolution of $\pmb{n}$ in time as expressed in the body frame can be written as follows: 
\begin{equation}
	\pmb{\dot{n}} = -\pmb{\omega}_B\times\pmb{n}
\end{equation}
In equilibrium, $\pmb{n}$ is fixed, therefore $\dot{\pmb{n}}=0$. This means that the fixed axis $\pmb{n}$ must be parallel to the angular velocity vector of the vehicle. In hovering state, assuming $\bar{\pmb{n}}$ is a unit vector ($|\bar{\pmb{n}}| = 1$), one can write the following (an overbar is used to denote equilibrium values): 
\begin{equation}
\label{nbar_eq}
    \bar{\pmb{n}} = \frac{\bar{\pmb{\omega}}_B}{|\bar{\pmb{\omega}}_B|}
\end{equation}
Total thrust force of the propeller and the streamline-shape fuselage help balance the vehicle's weight to keep it at constant altitude. This adds the following constraint to the system: 
\begin{equation}
\label{mg_const}
	|(\bar{\pmb{f}}_p+\bar{\pmb{f}}_B)\cdot\bar{\pmb{n}}| = m|\pmb{g}|
\end{equation}
Finally, by setting angular accelerations to zero, using two equations in (\ref{motor_eq}), three equations in (\ref{rot_eq}), three equations in (\ref{nbar_eq}) and one equation in (\ref{mg_const}), one obtains nine algebraic equations that can be solved for nine unknowns, which are: $\bar{i}$, $\bar{V}_m$, $\bar{\omega}_p$, $\bar{p}$, $\bar{q}$, $\bar{r}$, $\bar{n}_x$, $\bar{n}_y$ and $\bar{n}_z$. 

\section{OPTIMAL HOVER SOLUTIONS}

In this section, hover solution for various configurations of a monospinner is presented. Various attempts have been made to reduce power consumption in unmanned systems in general. Particularly, significant attempts in designing optimal control systems considering different constraints (e.g., on outputs, inputs, communication, etc) have been made \cite{khoji1}-\cite{bou}. However, the focus of this paper is on the mechanical/aerodynamic design and the objective is to find optimal-power hover solution for two specific configurations mentioned earlier.
First, we define a set of six design variables as follows: 
\begin{equation}\label{eq:x}
	x = (\alpha_p, \alpha_B, \frac{c_B}{c_p}, \frac{R_B}{R_p}, \delta, \frac{l}{R_B})
\end{equation}
where the first two variables are angles of attack for the propeller and fuselage, $\delta$ is the tilting angle of the rotor and the other three variables are related to the geometry of the vehicle. 

In order to meaningfully compare different configurations in terms of power consumption, we attempt to minimize the specific power (power-to-weight ratio) in hovering as follows: 
\begin{equation}
	\label{objective}
	P_s = \frac{\bar{P}_m}{m_T|\pmb{g}|} = f(x)
\end{equation}

The goal is to find the best set of parameters in $x$ that minimizes (\ref{objective}). 
\begin{equation}
\label{optim}
\begin{aligned}
& \underset{x}{\text{argmin}}
& & f(x) \\
& \text{s.t.} & &  equations \: of \: motion \\
\end{aligned}
\end{equation}

Given the complexity and nonlinearity of equations of motions, a numerical optimization routine (brute-force search) has been utilized, for which we used some specific physical parameters given below. One should note that the same physical parameters were used under the two aforementoned vehicle configurations.
\begin{equation}
\begin{matrix}
	c_p = 0.03 \: \text{m} \: \: \:,\: \: \:  R_p = 0.08  \: \text{m} \\
    g=9.81  \: \textrm{m/s$^2$} \: , \: \rho=1.225 \: \textrm{kg/m$^3$} \\
    K_{\tau_m}=0.02 \: \text{N.m/A} \: , \: K_v=0.02 \: \text{V/rad/s} \\
    \gamma = 9.75\times10^{-6}  \: \text{N.m.s/rad}
\end{matrix}
\end{equation}

\begin{figure}[t]
\centering
\includegraphics[scale=0.36]{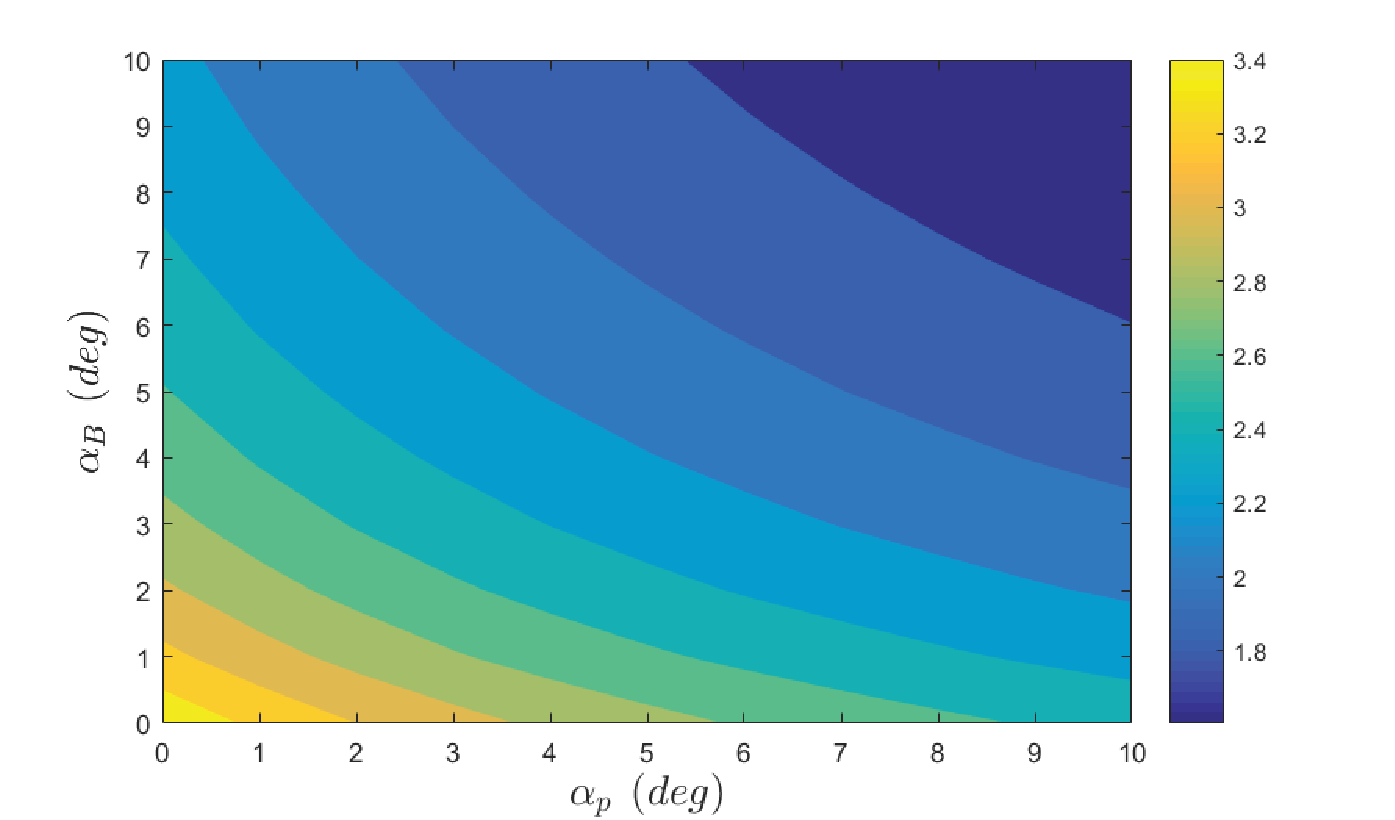} 
\caption{ Contours of specific power (units are in W/N) vs angle of attack of the fuselage and the propeller are presented. Blue represents low power and yellow represents high power. }
\label{coaxial1}
\end{figure}
\begin{figure}[t]
\centering
\includegraphics[scale=0.35]{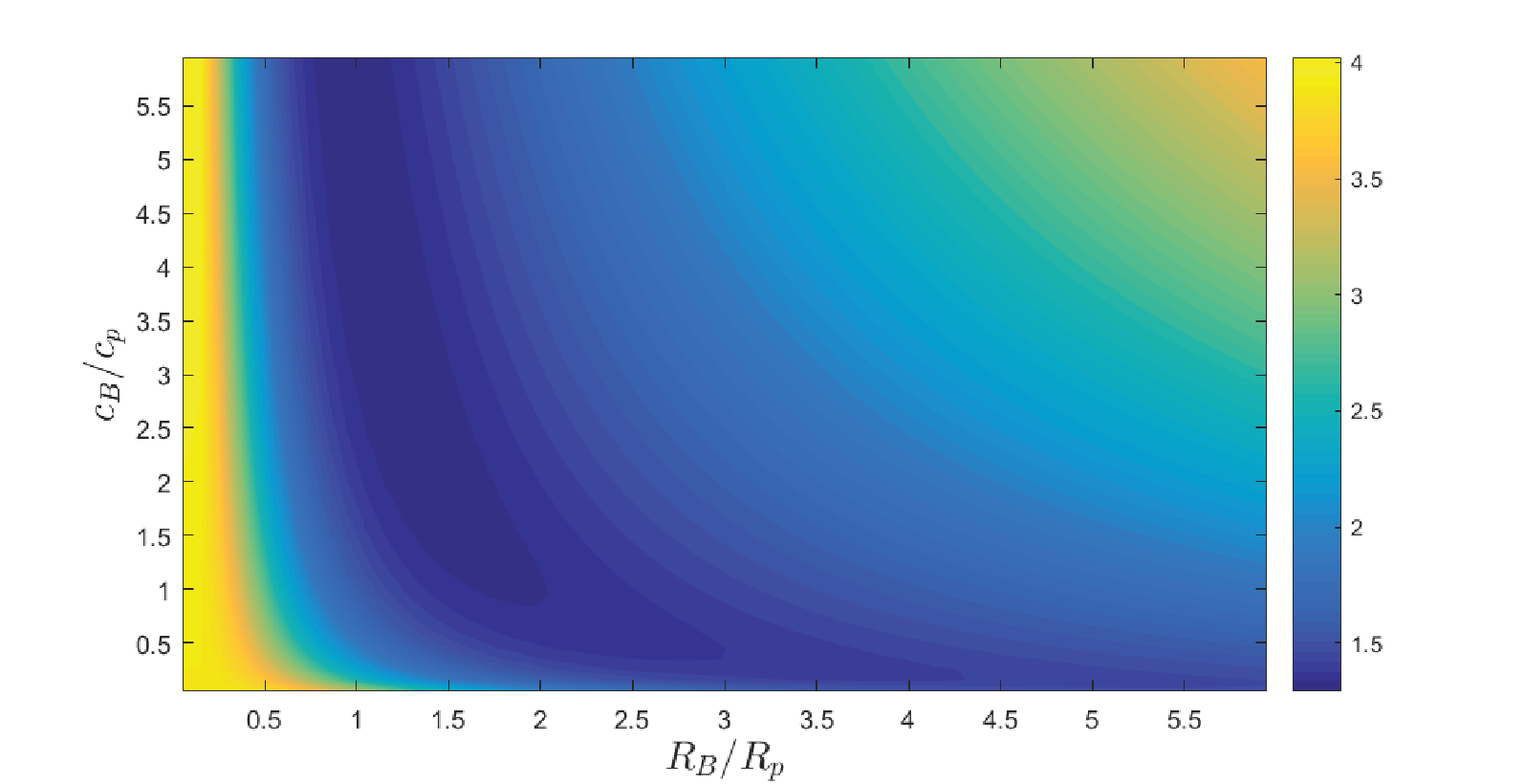} 
\caption{ Contours of specific power (units are in W/N) vs $c_B/c_p$ and $R_B/R_p$. Blue represents low power and yellow represents high power.}
\label{coaxial2}
\end{figure}

\subsection{First Configuration: Co-axial} 

In the co-axial configuration, the propeller's COM coincides with that in the fuselage ($l=0$) and the tilting angle $\delta$ is set to zero. It is assumed that the fuselage is streamline-shape and that the propeller and fuselage have a common axis of rotation. Unlike co-axial helicopters, in this configuration, the unbalanced moments in the system will turn the fuselage. In co-axial helicopters, some of the moments (e.g., reaction moment of the propeller) are used to balance the total moments in the system to bring angular velocities down to zero. In spinners, this requirement is relaxed, however, the angular velocity must remain constant instead of zero. Hover solution for different values of~\eqref{eq:x} is found and numerical results are presented in this section. 

Variations of the specific power versus angle of attack of the propeller and the fuselage are shown in Fig. \ref{coaxial1}. 

As shown in Fig. \ref{coaxial1}, for optimal utilization of the power, both fuselage and propeller should have the largest angle of attack (dark blue area). Also, from Fig.~\ref{coaxial2}, a set of $\frac{c_B}{c_p}$ and $\frac{R_B}{R_p}$ could be found such that the specific power of the co-axial spinner is minimized during hover (dark blue area). 

Note that under this co-axial configuration, because roll and pitch rates of the vehicle cannot be controlled, the x and y components of the position of the vehicle cannot be controlled. Therefore, this configuration would be capable of maintaining the altitude only. However, it can be shown that by tilting the rotor about the y-axis of the body frame, position control can be achieved as well. 

Because of the low dimension of the state space in \ref{objective}, an exhaustive search algorithm was used to find the minimum specific power hover solution. The resulting configuration is as follows: 
\begin{equation}
\label{coaxial_opt}
\begin{matrix}
	\frac{c_B}{c_p} = 1.05 \:\:,\:\: \frac{R_B}{R_p} = 1.75 \:\:,\:\: P_s = 1.3296 \: \text{W/N} \\
    \alpha_p = \alpha_B = \ang{10} \:\:,\:\: \frac{l}{R_B}=\delta=0
\end{matrix}
\end{equation}
And the resulting hover solution is found as follows: 
\begin{equation}
\begin{matrix}
	P_s = 1.3296 \: \text{W/N} \:\:,\:\: \bar{\omega}_p = 471.48 \: \text{rad/s}\\
    i_m = 0.25 \: \text{A} \:\:,\:\: \bar{\pmb{n}} = (0,0,1)^T \\
    \bar{\pmb{\omega}}_B = (0,0,-104.52)^T \: \text{rad/s} \:\:,\:\: V_m = 9.68 \: \text{V}
\end{matrix}
\end{equation}

This means that the optimal-power hover solution for the co-axial configuration can be reached when the fuselage is streamline-shape and aerodynamic. Without considering the streamline-shape fuselage, the optimal-power hover solution will be as follows: 
\begin{equation}
\begin{matrix}
	\frac{c_B}{c_p} = 0 \:\:,\:\: \frac{R_B}{R_p} = 0 \:\:,\:\: P_s = 3.926 \: \text{W/N} \\
    \alpha_p = \alpha_B = \ang{10} \:\:,\:\: \frac{l}{R_B}=\delta=0
\end{matrix}
\end{equation}

\begin{figure}[t]
\centering
\includegraphics[scale=0.39]{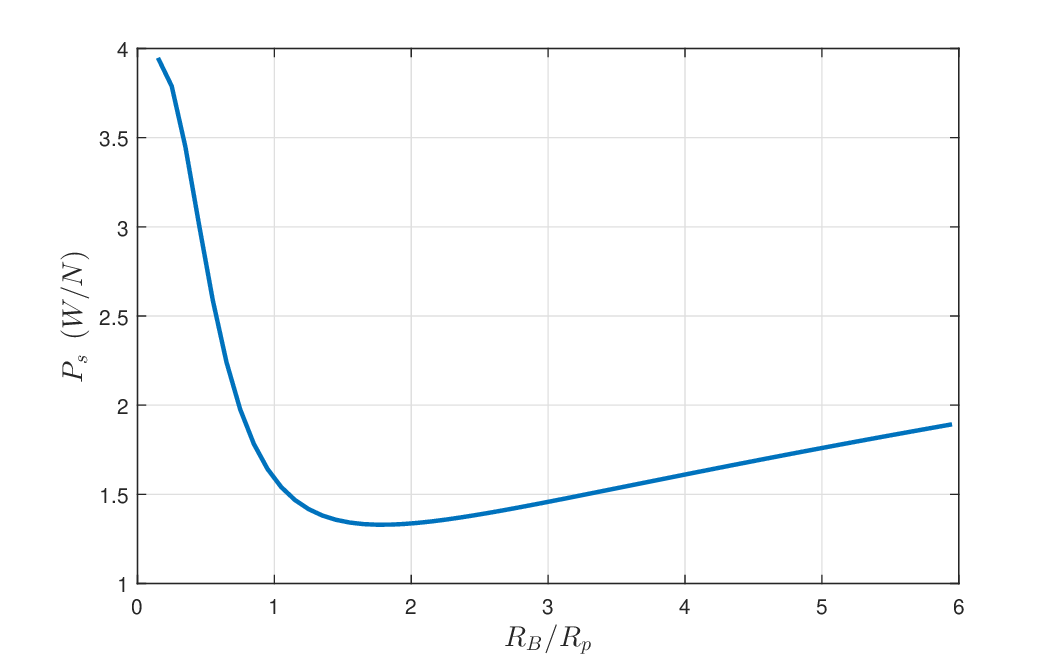} 
\caption{ Specific power vs $R_B/R_p$ when $c_B/c_p=1.05$}
\label{coaxial3}
\end{figure}

The results clearly show the significant contribution of the streamline-shape fuselage on reducing power consumption. Fig. \ref{coaxial3} also shows the strong relationship between the specific power and $\frac{R_B}{R_p}$ ratio.

\subsection{Second Configuration: Propeller is located at distance \textit{l} from COM of the Vehicle}

As a second alternative to the vehicle configuration, the propeller can be placed anywhere on the fuselage along the y-axis of the body frame. It is assumed that $l$ is a fraction of $R_B$ and that it lies on the y-axis of the body frame (could be either positive or negative). Like that in the previous configuration, it is assumed that the fuselage is streamline-shape and aerodynamically designed. Using the same design variables and solving the same optimization problem as in \ref{objective} and \ref{optim}, optimal-power hover solution was found.  

\begin{figure}[t]
\centering
\includegraphics[scale=0.39]{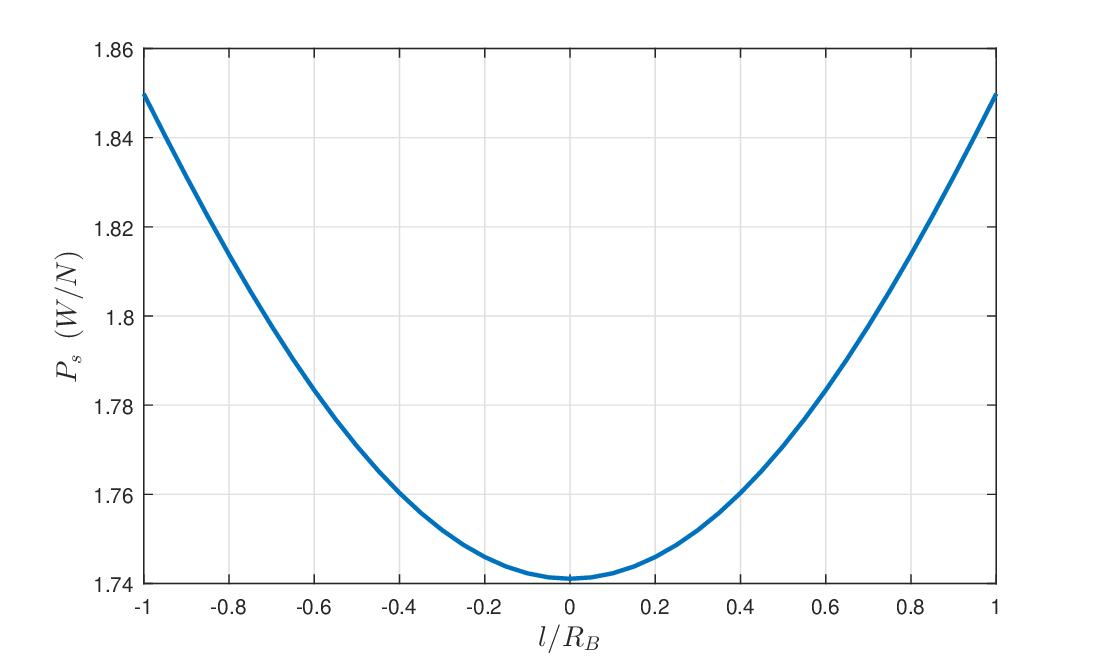} 
\caption{ Specific power vs $l$ for $\frac{c_B}{c_p}=1$ and $\frac{R_B}{R_p}=5$.}
\label{tilted_1}
\end{figure}

As an example, without considering the effects of tilting the rotor, specific power for hovering for all values of $l$ is presented in Fig. \ref{tilted_1}. It can be seen that the graph is symmetric, which means that the sign of $l$ does not affect the power consumption. Also, it can be conjectured that the propeller's COM should be aligen with that in the fuselage to obtain the minimum-power hover solution. 

Consequently, when $l$ and tilting angle $\delta$ are zero, the unit vector $\pmb{n}$ will be parallel to the gravitational acceleration. However, when $l$ is nonzero, only a component of $\pmb{n}$ will be parallel to gravitational acceleration and as a result, we should expect to have higher power consumption. 

Next, we consider the effect of tilting in the rotor. Incorporating the titling motion in the propellers would make the system highly nonlinear, which leads to the fact that, in some special cases, a feasible hover solution may not even exist. This also implies that the optimal solution would be very sensitive to the values of~\eqref{eq:x} at each iteration of the optimization.  
Here we present a numerical, but realistic, example, around the following initial configuration:  
\begin{equation}
	x_0 = (10,10,1,5,0,0),
\end{equation}
and solving the optimization problem to find the optimal-power hover solution results in the following configuration: 
\begin{equation}
	x = (10,10,0.95,5.19,0.16,1)
\end{equation}
and the hover solution for this configuration is as follows: 
\begin{equation}
\begin{matrix}
	P_s = 0.1325 \: \text{W/N}  \:\:,\:\: \bar{\omega}_p = -86.87 \: \text{rad/s}\\
    i_m = 0.17 \: \text{A} \:\:,\:\: \bar{\pmb{n}} = (0.002,0.017,0.999)^T \\
    \bar{\pmb{\omega}}_B = (-0.1,-0.66,-38.77)^T \: \text{rad/s} \:\:,\:\: V_m = 1.91 \: \text{V}
\end{matrix}
\end{equation}
It can be seen that the power consumption is reduced significantly. Because of the tilting angle, one component of the thrust force of the propeller generates fast yaw motion which consequently generates lift on the streamline-shape fuselage,  and it helps hover more efficiently. However, for large tilting angles, and due to nonlinearities, the numerical solver may not converge to an optimal solution. Currently, we are exploring about nonlinear solvers for this purpose such as the one in PETSc \cite{pouria}. 

\section{Discussion}

In this section a summary of the different configurations and effects of three factors on the performance and power consumption of the proposed spinner is presented. The outcomes for the first configuration (symmetric) can be summarized as: 
\begin{itemize}
\item By making the fuselage streamline-shape and aerodynamically  sound, one can use the fast rotations in the system to generate extra lift that leads to lower power consumption significantly. 
\item One should note that, without the tilting in the propellers, trajectory tracking cannot be achieved. 
\item The best chord and radius ratios between the propeller and the fuselage airfoils was found to obtain minimum-power hover solution.
\end{itemize}
For the second configuration (asymmetric), the outcomes can be summarized as follows: 
\begin{itemize}
\item The effect of moving the propeller's COM away from the fuselage's along the y-axis of the body frame was investigated and results show that the farther the distance, the higher power the vehicle will consume. 
\item The effects of the tilting angle on power consumption were presented. Results show that the tilting angle adds more nonlinearity in the system which makes the optimization solver sensitive to the initial guess.
\item It is shown that for small tilting angles, however, the power consumption is reduced. This was supported via one numerical example, as one representative of many cases tried in our lab. 
\end{itemize}

\section{CONCLUSIONS AND FUTURE WORK}

A comprehensive dynamic model of a mono-spinner UAV was presented. Two configurations were studied, namely (1)
symmetric, and (2) asymmetric. The effect of three design factors on the overall power consumption and flight stability were studied that include: (1) the aerodynamic shape of the fuselage, (2) the distance between the rotor and the fuselage's COM, and (3) the titling angle between the rotor and the fuselage's principle axis.

Hover solutions for the two configurations were presented and Compared. The concept of specific power
was used via numerical optimizations for comparison purposes. Results show that power-optimal solutions exist. However, while a global optimal solution is achievable under the first configuration, a number of locally-optimal solutions would exist under the second configuration due to non-convexity of the optimization objective function. Also results show that using streamline-shape aerodynamic fuselage has significant effects on power-consumption. This paper only presents the results under the assumption of a hovering flight at steady state. In addition to experimental validation, future plans include: investigating about
the power consumption when going through a flight state transition via conventional maneuvers such as
tracking trajectories, servoing towards an object, point to point motion, or landing on a desired spot.

\addtolength{\textheight}{-12cm}   % This command serves to balance the column lengths

%%%%%%%%%%%%%%%%%%%%%%%%%%%%%%%%%%%%%%%%%%%%%%%%%%%%%%%%%%%%%%%%%%%%%%%%%%%%%%%%

%%%%%%%%%%%%%%%%%%%%%%%%%%%%%%%%%%%%%%%%%%%%%%%%%%%%%%%%%%%%%%%%%%%%%%%%%%%%%%%%

%%%%%%%%%%%%%%%%%%%%%%%%%%%%%%%%%%%%%%%%%%%%%%%%%%%%%%%%%%%%%%%%%%%%%%%%%%%%%%%%

%%%%%%%%%%%%%%%%%%%%%%%%%%%%%%%%%%%%%%%%%%%%%%%%%%%%%%%%%%%%%%%%%%%%%%%%%%%%%%%%

\end{document}